\documentclass[11pt,dvipdfmx]{amsart}
\usepackage{txfonts,mathptmx}
\usepackage{mathrsfs}
\usepackage[dvipdfmx,colorlinks=true]{hyperref}
\usepackage[msc-links]{amsrefs}

\usepackage{tikz-cd}

\numberwithin{equation}{section}

\newtheorem*{thm-nn}{Theorem}

\newtheorem*{cor-nn}{Corollary}

\theoremstyle{definition}

\theoremstyle{remark}

\newtheorem*{ack}{Acknowledgement}

\newcommand{\qdl}{\ast}

\newcommand{\Z}{\mathbb{Z}}

\newcommand{\map}{\rightarrow}

\makeatletter
\@namedef{subjclassname@2020}{%
  \textup{2020} Mathematics Subject Classification}
\makeatother

\BibSpec{arXiv}{%
  +{}{\PrintAuthors}{author}
  +{,}{ \textit}{title}
  +{}{ \parenthesize}{date}
  +{,}{ arXiv }{eprint}
  +{,}{ primary class }{primaryclass}
}

\begin{document}
\title[Embedding Alexander quandles]{
Embedding Alexander quandles into groups}
\author[T. Akita]{Toshiyuki Akita}
\address{Department of Mathematics, Faculty of Science, Hokkaido University,
Sapporo, 060-0810 Japan}
\email{akita@math.sci.hokudai.ac.jp}




\maketitle

\begin{abstract}
For any twisted conjugate quandle $Q$, 
and in particular any Alexander quandle,
there exists a group $G$ such that $Q$ is embedded
into the conjugation quandle of $G$.
\end{abstract}

\section{Embeddable quandles}
A non-empty set $Q$ equipped with a binary operation $Q\times Q\map Q$,
$(x,y)\mapsto x\qdl y$ is called a \emph{quandle}
if it satisfies the following three axioms:
\begin{enumerate}
\item $x\qdl x=x$ $(x\in Q)$,
\item $(x\qdl y)\qdl z=(x\qdl z)\qdl (y\qdl z)$ $(x,y,z\in Q)$,
\item For all $x\in Q$,
the map $S_x\colon Q\map Q$ defined by $y\mapsto y\qdl x$ is bijective.
\end{enumerate}
Quandles were introduced 
independently by Joyce \cite{MR638121} and Matveev \cite{MR672410}.
Since then, quandles have been important objects in the study of 
knots and links,
set-theoretical solutions of the Yang-Baxter equation,
Hopf algebras and many others.
We refer to Nosaka \cite{nosaka-book} for further details of quandles.

A map $f \colon Q \rightarrow Q'$ of quandles is called a
\emph{quandle homomorphism} if it satisfies
$f(x\qdl y)=f(x)\qdl f(y)$ $(x,y\in Q)$.
Given a group $G$, the set $G$ 
 equipped with a quandle operation
$h\qdl g\coloneqq g^{-1}hg$ is called the \emph{conjugation quandle}
of $G$ 
and is denoted by $\mathrm{Conj}(G)$.
A quandle $Q$ is called \emph{embeddable} if
there exists a group $G$ and an injective quandle homomorphism
$Q\to\mathrm{Conj}(G)$.
Not all quandles are embeddable (see the bottom of \S\ref{sec:tcq}).

In their paper \cite{MR3718201},
Bardakov-Dey-Singh proposed the question
``For which quandles $X$ does there exists a group $G$
such that $X$ embeds in the conjugation quandle $\mathrm{Conj}(G)$?''
\cite{MR3718201}*{Question 3.1},
and proved that
Alexander quandles associated with fixed-point free involutions
are embeddable
\cite{MR3718201}*{Proposition 3.2}.
The following is a list of embeddable quandles of which the author is aware:
(1) free quandles and free $n$-quandles (Joyce 
\cite{MR638121}*{Theorem 4.1 and Corollary 10.3}),
(2) commutative quandles, latin quandles and simple quandles
(Bardakov-Nasybullov \cite{MR4129183}*{\S5}),
(3) core quandles (Bergman \cite{MR4255023}*{(6.5)}),
(4) generalized Alexander quandles associated with fixed-point free
automorphisms (Dhanwani-Raundal-Singh \cite{dehn}*{Proposition 3.12}),
and (5) free \emph{c}-nilpotent quandles (Darn\'e
\cite{nilpotent}*{Proposition 2.18}).

In this short note, we will show that twisted conjugation quandles,
which include all Alexander quandles, are embeddable, thereby
generalize the aforementioned result of Bardakov-Dey-Singh.

\section{Embeddings of twisted conjugation quandles}\label{sec:tcq}
Let $G$ be an additive abelian group and let $\phi\colon G\to G$ be a
group automorphism of $G$.
The \emph{Alexander quandle} $\mathrm{Alex}(G,\phi)$
associated with $\phi$ is the
set $G$ equipped with the quandle operation
\[
g*h\coloneqq\phi(g)+h-\phi(h).
\]
Let $G$ be a group and let $\phi\colon G\to G$ be an
automorphism of $G$.
The \emph{twisted conjugation quandle} $\mathrm{Conj}(G,\phi)$
associated with $\phi$ is the the set $G$ 
equipped with the quandle operation
\[
g\ast h\coloneqq \phi(h^{-1}g)h.
\]
Observe that an Alexander quandle $\mathrm{Alex}(G,\phi)$
is precisely a twisted conjugation quandle $\mathrm{Conj}(G,\phi)$
whose underlying group $G$  is abelian.
Twisted conjugation quandles appeared in
Andruskiewitsch-Gra\~{n}a \cite{MR1994219}*{\S1.3.7}
under the name \emph{twisted homogeneous crossed sets}.
We prefer the name twisted conjugation quandles
because $\mathrm{Conj}(G,\phi)=\mathrm{Conj}(G)$
if $\phi$ is the identity map.
It should be emphasized that $\mathrm{Conj}(G,\phi)$ is
different from the \emph{generalized Alexander quandle}
associated with $(G,\phi)$.
The latter has the same underlying set $G$, but 
with the different quandle operation 
$g\ast h\coloneqq \phi(gh^{-1})h$.
Now we prove that $\mathrm{Conj}(G,\phi)$ is embeddable:
\begin{thm-nn}
Any twisted conjugation quandle 
is embeddable. In particular, any Alexander quandle
is embeddable.
\end{thm-nn}
\begin{proof}
Given a twisted conjugation quandle $\mathrm{Conj}(G,\phi)$,
we will construct
an explicit embedding $\mathrm{Conj}(G,\phi)
\to\mathrm{Conj}(H)$.
Let $\Z$ be the additive group of integers, and
let $H\coloneqq G\rtimes_{\phi}\Z$ be the semidirect product of
$G$ and $\Z$ associated with $\phi$.
Namely, $H$ equals to $G\times\Z$ as sets.
The group law on $H$ is given by
\[ (g,m)\cdot (h,n)\coloneqq (\phi^{n}(g)h,m+n).\]
The inverse of $(g,m)\in H$ is 
\[(g,m)^{-1}=(\phi^{-m}(g^{-1}),-m).\]
Observe that
\begin{align*}
(g,1)\ast (h,1)&\coloneqq (h,1)^{-1}\cdot (g,1)\cdot (h,1)
=(\phi^{-1}(h^{-1}),-1)\cdot (g,1)\cdot (h,1)\\
&=(\phi^{-1}(h^{-1}),-1)\cdot (\phi(g)h,2)
=(\phi^{2}(\phi^{-1}(h^{-1}))\phi(g)h,1)\\
&=(\phi(h^{-1})\phi(g)h,1)=(\phi(h^{-1}g)h,1)
\end{align*}
holds in $\mathrm{Conj}(H)$, and
we conclude that the injective map $G\to H$ defined by $g\mapsto (g,1)$
is an injective quandle homomorphism $\mathrm{Conj}(G,\phi)
\to\mathrm{Conj}(H)$, hence verifying the theorem.
\end{proof}
Now let $Q$ be an arbitrary quandle. 
The \emph{associated group} $\mathrm{As}(Q)$ of $Q$ is
the group defined by the presentation
\[
\mathrm{As}(Q)\coloneqq\langle e_{x}\,(x\in Q)\mid
e_{y}^{-1}e_{x}e_{y}=e_{x\ast y}\,(x,y\in Q)\rangle.
\]
A quandle $Q$ is called \emph{injective} if the canonical map
$Q\to\mathrm{As}(Q)$ defined by $x\mapsto e_{x}$
$(x\in Q)$ is injective.
The injectivity of finite quandles is important in the study
of set-theoretical solutions of the Yang-Baxter equation
(see Lebed-Vendramin \cite{MR3974961} for instance).
According to Joyce \cite{MR638121}*{Section 6}
(see also Dhanwani-Raundal-Singh \cite{dehn}*{Theorem 3.8}),
a quandle $Q$
is injective if and only if $Q$ is embeddable.
As a byproduct of the theorem, we obtain the following
corollary:
\begin{cor-nn}
Any twisted conjugation quandle 
is injective. In particular, any Alexander quandle
is injective. 
\end{cor-nn}
Finally, we remark that not all quandles are embeddable.
Indeed, there exist quandles which are not injective 
and hence are not embeddable.
See Joyce \cite{MR638121}*{Section 6} and 
Bardakov-Nasybullov \cite{MR4129183}*{\S4} for
examples of such quandles.

\begin{ack}
The author would like to thank the anonymous referee for valuable
comments improving this paper.
The author was partially supported by JSPS KAKENHI Grant Number
20K03600.
\end{ack}

\end{document}